\documentclass{article}
\usepackage{amssymb}

\newtheorem{satz}{Satz}[section]
\newtheorem{corollary}[satz]{Corollary}
\newtheorem{example}[satz]{Example}
\newtheorem{proposition}[satz]{Proposition}
\newtheorem{remark}[satz]{Remark}
\newenvironment{varthm*}[1]{\begin{list}{}{\labelwidth=0cm \leftmargin=0cm}\item[\hspace{\labelsep}\bfseries #1.]\itshape}{\end{list}}

\newcommand\qedsymbol{\frame{\rule[0pt]{0pt}{8pt}\rule[0pt]{8pt}{0pt}}}
\newcommand\proofname{Proof}
\newenvironment{proof}[1][\proofname]%
   {\trivlist\item[\hskip\labelsep{\it #1.}]}
   {\hspace*{\fill}\qedsymbol\endtrivlist}

\newcommand\subjclass[1]{{\renewcommand\thefootnote{}\footnotetext{2000 \textit{Mathematics Subject Classification:} #1.}}}
\newcommand\keywords[1]{{\renewcommand\thefootnote{}\footnotetext{\textit{Keywords.} #1.}}}

\renewcommand\ge{\geqslant}  
\renewcommand\le{\leqslant}  
\renewcommand\O{\mathcal O}
\newcommand\restr[1]{\big|_{#1}}
\newcommand\set[1]{\left\{\,#1\,\right\}}
\newcommand\with{\ \vrule\ }
\newcommand\midtext[1]{\quad\mbox{#1}\quad}
\newcommand\be{\begin{eqnarray*}}
\newcommand\ee{\end{eqnarray*}}
\newcommand\eqnref[1]{(\ref{#1})}
\newcommand\eqdef{=_{\rm def}}
\newcommand\eps{\varepsilon}
\newcommand\Q{\mathbb Q}
\newcommand\R{\mathbb R}
\newcommand\Z{\mathbb Z}
\newcommand\engqq[1]{``#1''}
\newcommand\tensor{\otimes}
\newcommand\inverse{^{-1}}
\newcommand\Biginverse{^{\hskip-0.25em-1}}
\renewenvironment{cases}{\left\{\begin{array}{cl}}{\end{array}\right.}
\newcommand\binom[2]{{#1\choose#2}}
\newcommand\longto{\longrightarrow}

\newcommand\newop[2]{\newcommand#1{\mathop{\rm #2}\nolimits}}
\newop\mult{mult} 
\newop\NS{NS} 
\newop\Bl{Bl} 
\newop\Num{Num}
\newop\End{End}
\newop\Amp{Amp}
\newop\dAmp{\partial\!\Amp}
\newop\id{id}
\def\slope(#1,#2){\sigma(#1,#2)}

\newcommand\J{\mathcal J}

\begin{document}

\title{A criterion for an abelian variety to be simple}
\author{Thomas Bauer\thanks{The author was supported by DFG grant BA 1559/4-3.}}
\date{October 16, 2007}
\maketitle

\begin{abstract}

   In this note we give a numerical 
   criterion that expresses the condition that an abelian variety
   be simple in terms of an invariant that is closely related 
   to the $s$-invariant of Ein-Cutkosky-Lazarsfeld.
\end{abstract}

\section*{Introduction}

   An abelian variety is \textit{simple} if it does not contain
   any non-trivial abelian subvarieties. 
   The purpose of this note is to 
   provide a
   criterion 
   that expresses simpleness 
   as a condition on
   the codimension one level:

\begin{varthm*}{Theorem}
   Let $(X,L)$ be a polarized 
   abelian variety over the complex numbers. 
   The following statements are equivalent:

\begin{itemize}
\item[\rm(i)]
   For every line bundle
   $M$ on $X$ that is not proportional to $L$ 
   in $\Num(X)\tensor\Q$, the supremum
   $$
      \sup\set{t\in\R\with L-tM\mbox{ is nef}}
   $$
   is irrational or equal to $\infty$.
\item[\rm(ii)]
   $X$ is simple.
\end{itemize}
\end{varthm*}
   Here $\Num(X)$ denotes the group of numerical
   equivalence classes of line bundles on~$X$.
   A line bundle $B$ is \textit{nef} if $B\cdot C\ge 0$ for every
   irreducible curve $C$ in $X$.

   The criterion provides new examples where the 
   $s$-invariants of Cutkosky-Ein-Lazarsfeld are irrational 
   (see Example~\ref{example}).
   It may also be viewed as a statement about the geometry of the
   ample cone of $X$ (see Remark~\ref{remark}).

\paragraph*{\it Conventions.}
   We work throughout over the field of
   complex numbers. 
   Additive notation will be used for the tensor product of line
   bundles.

\section{Proof of the theorem}

   Let $X$ be an abelian variety. For line bundles $L$ and $M$ 
   on $X$ we consider the number
   \begin{equation}\label{slope}
      \slope(L,M)\eqdef\sup\set{t\in\R\with L-tM\mbox{ is nef}}
      \in\R\cup\{\pm\infty\} 
   \end{equation}
   We start by giving an algebraic characterization of
   $\slope(L,M)$, when $L$ is ample.

\begin{proposition}\label{crit s}
   Let $(X,L)$ be a polarized abelian variety of dimension $n$,
   and let $M$ be any line
   bundle on $X$.
   If $\zeta(L,M)$ denotes the maximal root of the polynomial
   $$
      \chi(uL-M)=\frac{1}{n!}(uL-M)^n\in\Q[u]
   $$
   (all of whose roots are in any event real), then
   $$
      \slope(L,M)=\begin{cases}
                \infty & \mbox{if } \zeta(L,M)\le 0, \\[\medskipamount]
                \displaystyle\frac1{\zeta(L,M)} & \mbox{if } \zeta(L,M)>0. \
      \end{cases}
   $$
\end{proposition}

\proof
   We will make use of the isomorphism
   of $\Q$-vector spaces
   \be
     \NS_\Q(X) &\longto& \End^s_\Q(X) \\
     B &\longmapsto&f_B:=\phi\inverse_L\circ\phi_B \ ,
   \ee
   where $\phi_L$ and $\phi_B$ are the canonical homomorphisms
   $X\to\widehat X$ to the dual abelian variety $\widehat X$.
   Note first that $B$ is a nef class if and only if the $\Q$-endomorphism
   $f_B$ (or, more precisely, the derivative $d(mf_B):T_0X\to T_0X$
   of a suitable multiple $mf_B$, $m>0$)
   has no negative eigenvalues.
   In fact, the characteristic polynomial $P_B$ of $f_B$ satisfies
   $$
      P_B(m)=\frac{\chi(mL-B)}{\chi(L)}
   $$
   for $m\in\Z$ (see \cite[Sect.\ 5.2]{LB}), so that
   its alternating coefficients
   are positive multiples of the intersection numbers
   $$
      L^kB^{n-k} \qquad(0\le k\le n) \ ,
   $$
   where $n=\dim X$. But $B$ is nef if and only if all $L^kB^{n-k}$
   are non-negative (see e.g.\ \cite[Lemma 1.1]{Bau98}), and
   -- since $f_B$ is symmetric and has therefore only real eigenvalues --
   this is equivalent to saying that $P_B$ has no negative roots.

   So in particular, for $t\in\Q$ the line bundle $L-tM$ is nef if and
   only if the endomorphism
   $$
      f_{L-tM}=\id_X-tf_M
   $$
   has no negative eigenvalues.  But the eigenvalues of $\id_X-tf_M$
   are the numbers $1-t\lambda$ where $\lambda$ runs through the set of
   eigenvalues of $f_M$.
   From this the assertion follows.
\endproof

\begin{proof}[Proof of the theorem]
   Suppose first that (i) holds and assume by way of
   contradiction that there is
   an abelian subvariety
   $Y\subset X$ different from $X$ and $0$.
   Consider the norm endomorphism
   $N_Y\in\End(X)$
   of $Y$ with
   respect to the polarization 
   $L$ (see \cite{BirLan91} or \cite[Sect.\ 5.3]{LB} for details
   on norm endomorphisms).  The pullback
   $$
      M \eqdef N_Y^*L
   $$
   then corresponds to the endomorphism $f_M=N_Y^2=eN_Y$, where
   $e$ is the exponent of the induced polarization $L\restr Y$,
   i.e., the minimal positive integer $e$
   such that
   $$
      e\phi\inverse_{L\restr Y} \in\End(X)
   $$
   is an (integral) homomorphism.
   Since $L$ is ample and $M$ is non-trivial and nef, but certainly not ample,
   the bundles $L$ and $M$ are not proportional in
   $\Num(X)\tensor\Q$.
   Now, the eigenvalues of $N_Y$ are 0 and $e$, and hence
   Proposition~\ref{crit s} implies that
   $$
      \slope(L,M)=\frac1{e^2}\in\Q \ ,
   $$
   which is a contradiction with (i).

   Supposing now that (i) does not hold, we will show that 
   (ii) does not hold as well.
   We can 
   argue
   as in the first part of the
   proof of \cite[Proposition 1.2]{Bau98}.  In brief,
   if $\slope(L,M)$ is rational, then
   a suitable rational
   multiple of the class $L-\slope(L,M)M$ is an integral class $B$
   on $X$, which is nef but not ample.  The homomorphism $\phi_B$ has
   therefore a non-trivial kernel.
   On the other hand, since $L$ and $M$ are not proportional,
   $B$ cannot be topologically
   trivial, and consequently $\phi_B$ cannot be the zero morphism.
   So the connected component of its
   kernel containing the point~0 is a non-trivial abelian subvariety
   of $X$.
   This completes the proof.
\end{proof}

\section{Complements and application to $s$-invariants}

   We give here two further 
   applications of Proposition \ref{crit s}, and we point out the
   geometric consequences of the theorem.
   Our first observation says in effect that if
   $\slope(L,M)$ is rational, then there are only finitely many
   possibilities for its value. 
   
\begin{corollary}\label{div cond}
   Let $(X,L)$ be a polarized abelian variety of dimension $n$,
   and let $M$ be a line
   bundle on $X$ such that $-M$ is not nef.
   If $\slope(L,M)$ is a rational number, then it is of the form
   $$
      \slope(L,M)=\frac pq \ ,
   $$
   where $p$ and $q$ are coprime integers satisfying the divisibility conditions
   $$
      p|L^n \midtext{and} q|M^n \ .
   $$
\end{corollary}

\begin{proof}
   Write $\slope(L,M)=p/q$ with coprime integers $p$ and $q$.
   By Proposition~\ref{crit s}, the rational number $q/p$ is the
   maximal root of the polynomial $\chi(uL-M)$.  Now, the polynomial
   $$
      n!\cdot\chi(uL-M)
   $$
   has integer coefficients, its leading coefficient
   is $L^n$, and the constant term is (up to a possible sign) $M^n$.
   This implies the assertion.
\end{proof}

\begin{example}[Irrational $s$-invariants]\label{example}\rm 
   We establish here the relationship between our result and the
   $s$-invariants introduced by Cutkosky-Ein-Lazarsfeld in
   \cite{CEL01}. In particular we obtain many new examples
   of irrational $s$-invariants. 
   
   Consider a coherent ideal sheaf $\J\subset\O_X$ on a smooth
   projective variety $X$, and let $\nu:Y\to X$ be the blow-up of
   $X$ along $\J$. We have $\J\cdot\O_Y=\O_Y(-F)$ for an
   effective Cartier divisor $F$ on $Y$. Fixing an ample divisor
   $H$ on $Y$, the \textit{$s$-invariant}
   of $\J$ with respect to $H$ is
   defined to be the positive real number
   $$
      s_H(\J)=\min\set{s\in\R\with s\cdot\nu^*H-F\mbox{ is nef}} 
   $$
   (see \cite[Sect.~1]{CEL01}).
   Interestingly, the $s$-invariant governs (among other things) 
   the asymptotic
   regularity of powers of $\J$ (see \cite[Sect.~3]{CEL01}).
   When $\J$ is the ideal sheaf of a point $x\in X$, the
   reciprocal of $s_H(\J)$ is the \textit{Seshadri constant}
   $\eps(H,x)$,
   as introduced by Demailly
   (see \cite{Dem92} and \cite[Chapt.~5]{PAG}). 
   Paoletti (\cite{Pao94}, \cite{Pao95}) has studied the case where 
   $\J$ is the ideal sheaf of a smooth curve in a threefold. 
   It is natural to ask whether $s$-invariants can become
   irrational. 
   While it is still unknown whether 
   this can happen for 
   Seshadri constants (i.e.\ when $\J$ is the ideal sheaf of a
   point), it does happen for $s$-invariants in general.
   The first examples, due to 
   Ein-Cutkosky-Lazarsfeld,
   are $s$-invariants of curves on
   suitable abelian surfaces (see \cite[Example~1.7 and
   Example~1.11]{CEL01}).
   Our result clarifies the picture on abelian varieties 
   in the following way.
   If $(X,L)$ is a polarized abelian variety and $D$ an effective
   divisor on $X$, then the number $\sigma(L,D)$ defined in
   \eqnref{slope} is just the
   reciprocal of an $s$-invariant,
   $$
      \sigma(L,D)=\frac1{s_L(\J_{D/X})} \ ,
   $$
   where $\J_{D/X}$ is 
   is the ideal sheaf of $D$ in $X$.
   So the present result implies that on simple abelian
   varieties in fact \textit{all} such $s$-invariants
   are irrational (as long as $L$ and $\O_X(D)$ are not numerically
   proportional), while on non-simple abelian varieties 
   rational $s$-invariants occur.
\end{example}

\begin{remark}\label{remark}\rm
   It may also be useful to think of the theorem -- and the
   invariant $\slope(L,M)$ -- in terms of the geometry of the
   ample cone as follows. Given a polarized abelian variety
   $(X,L)$, we may ask \engqq{how far} $L$ is away from the
   boundary $\dAmp(X)$ of the ample cone of $X$. The number
   $\sigma(L,M)$ is then just the distance of $L$ to $\dAmp(X)$
   in the direction of $-M$, measured in units of $M$. When $X$
   contains a non-trivial abelian subvariety $Y$ and $M=N_Y^*L$,
   then the proof of the theorem tells us that this distance is
   $\le 1$, and in fact a rational number of the form $1/e^2$.
\end{remark}

   Finally, we establish a lower bound 
   on $\slope(L,M)$ that can
   be computed explicitly from the intersection numbers $L^kM^{n-k}$.
   In view of Example \ref{example}, this 
   gives an upper bound on the corresponding $s$-invariant.

\begin{corollary}\label{estimate}
   Let $(X,L)$ be a polarized abelian variety of dimension $n$,
   and let $M$ be a line
   bundle on $X$ such that $-M$ is not nef. Then
   $$
        \slope(L,M)\ge\left(\displaystyle 1+\max\limits_{0\le k<n} \binom
        nk\frac{L^kM^{n-k}}{L^n}\right)\Biginverse
   $$
\end{corollary}

\begin{proof}
   By Riemann-Roch, the coefficient $a_k$ at $u^k$
   in the polynomial
   $\chi(uL-M)$ is given by
   $$
      a_k=(-1)^{n-k}\frac{L^kM^{n-k}}{k!(n-k)!}
   $$
   for $0\le k\le n$.
   It is a theorem of Cauchy (see for instance
   \cite[Theorem 27.2]{Mar66}) that
   all roots of a complex polynomial $P(u)=\sum_{k=0}^n a_ku^k$,
   with $a_n\ne 0$,
   lie within the circle around 0 of radius
   $$
      1+\max\limits_{0\le k<n} \left|\frac{a_k}{a_n}\right| \ ,
   $$
   so that in our case we have in particular
   $$
      \zeta(L,M)<1+\max\limits_{0\le k<n} \binom nk\frac{L^kM^{n-k}}{L^n}
   $$
   for the maximal root $\zeta(L,M)$ of $\chi(uL-M)$.
   Since $-M$ is not nef, $\zeta(L,M)$ is positive,
   and hence Proposition~\ref{crit s} implies the assertion.
\end{proof}

   When 
   Proposition~\ref{crit s} is being used
   to bound 
   $\slope(L,M)$, the issue is to estimate the roots of
   polynomials in terms of their coefficients. The estimate used 
   in the proof of Corollary~\ref{estimate}
   is among the most immediate 
   ones; the reader may consult for instance
   \cite{Mar66} for more refined estimates.

\bigskip\noindent
   Fach\-be\-reich Ma\-the\-ma\-tik und In\-for\-ma\-tik \\
   Philipps-Uni\-ver\-si\-t\"at Mar\-burg \\
   Hans-Meer\-wein-Stra{\ss}e \\
   D-35032~Mar\-burg, Germany \\
   E-Mail: tbauer@mathematik.uni-marburg.de

\end{document}